\documentclass[10pt, twocolumn]{article}

\usepackage{authblk}
\usepackage{mathrsfs}
\usepackage{amsmath}
\usepackage{amssymb}
\newcommand{\ol}{\overline}
\newcommand{\wt}{\widetilde}
\newcommand{\vocab}[1]{\textbf{#1}}

\usepackage{natbib}
\usepackage{amsthm}
\newtheorem{thm}{Theorem}
\newtheorem{ex}{Example}

\newtheorem{de}{Definition}
\newtheorem{rem}{Remark}

\makeatletter
\renewcommand{\maketitle}{\bgroup\setlength{\parindent}{0pt}
\begin{flushleft}
  \textbf{\@title}

  \@author
\end{flushleft}\egroup
}
\makeatother

\title{\Huge Reflective prolate-spheroidal operators and the KP/KdV equations}

\author[a]{\large W. Riley Casper}
\author[b]{F. Alberto Gr{\"u}nbaum}
\author[a]{Milen Yakimov}
\author[c]{Ignacio Zurri\'an}

\affil[a]{\normalsize Dept of Mathematics, Louisiana State University, Baton Rouge, LA 70803, U.S.A.}
\affil[b]{Dept of Mathematics, University of California, Berkeley, CA 94720, U.S.A.}
\affil[c]{FaMAF, Universidad Nacional de C\'ordoba, 5000 C\'ordoba, Argentina}

\date{\today}

\usepackage[showframe=false, margin=.60in]{geometry}

\usepackage{abstract}
\renewcommand{\abstractname}{}    % clear the title
\usepackage{lipsum}

\usepackage{afterpage}

\renewenvironment{abstract}
 {\small
  \begin{center}
  \bfseries \abstractname\vspace{-.5em}\vspace{0pt}
  \end{center}
  \list{}{%
    \setlength{\leftmargin}{10mm}% <---------- Change margin here
    \setlength{\rightmargin}{\leftmargin}%
  }%
  \item\relax}
 {\endlist}

\usepackage[utf8]{inputenc}

\usepackage{wrapfig}
\usepackage[lofdepth,lotdepth]{subfig}
\usepackage{booktabs}

\usepackage{rotating}

\usepackage{caption}

\usepackage{csquotes}

\usepackage{todonotes}

\usepackage{stfloats}

\usepackage[T1]{fontenc}
\usepackage{textcomp}
\usepackage{times}

\usepackage{framed} %um Boxen zu machen

\usepackage{setspace}
\usepackage{color}
\definecolor{black}{gray}{0} % 10% gray

\usepackage{tabularx}
\newcolumntype{b}{X}
\newcolumntype{s}{>{\hsize=.5\hsize}X}

\begin{document}

\twocolumn[
\begin{@twocolumnfalse}
\maketitle
\begin{abstract}
Commuting integral and differential operators connect the topics of Signal Processing, Random Matrix Theory, and Integrable Systems.
Previously, the construction of such pairs was based on direct calculation and concerned concrete special cases, leaving behind important families such as the 
operators associated to the rational solutions of the KdV equation. We prove a general theorem that the integral operator 
associated to every wave function in the infinite dimensional Adelic Grassmannian ${\mathrm{Gr}}^{\mathrm{ad}}$ of Wilson always reflects a differential operator (in the sense of Definition \ref{Def-reflect} below). 
This intrinsic property is shown to follow from the symmetries of Grassmannians of KP wave functions, 
where the direct commutativity property holds for operators associated to wave functions fixed by Wilson's sign involution but is violated in general.

Based on this result, we prove a second main theorem that the integral operators in the computation of the singular values of 
the truncated generalized Laplace transforms associated to all bispectral wave functions of rank 1 reflect a differential operator. 
A $90^\circ$ rotation argument is used to prove a third main theorem that the integral operators in the computation of the singular values of 
the truncated generalized Fourier transforms associated to all such KP wave functions commute with a differential operator. 
These methods produce vast collections of integral operators with prolate-spheroidal properties, including as special cases the 
integral operators associated to all rational solutions of the KdV and KP hierarchies considered by Airault-McKean-Moser and Krichever, respectively, in the late 70's. 
Many novel examples are presented.
\end{abstract}
%\keywords{Prolate-spheroidal integral operators $|$ Reflectivity $|$ Rational solutions of the KdV and KP equations}\newline
\vspace{0.1in}
\end{@twocolumnfalse}
]

\section{Background}
\subsection{Commuting integral and differential operators}
%\pnasdoi{\small \url{www.pnas.org/cgi/doi/10.1073/pnas.1906098116}}
In a pair of ground-breaking works from the late 1940's Claude Shannon laid down the mathematical foundations of communication theory \cite{Shannon1,Shannon2}. 
One of the key problems which he raised was: What is the best information that one can infer for a signal $f(t)$ which is time limited to the interval $[-\tau,\tau]$ 
from knowing its frequencies in the interval $[-\kappa,\kappa]$? This double concentration problem leads to the study of the singular values of an operator given by a finite Fourier transform
$$
(Ef)(z)=\int_{-\tau}^\tau e^{izx}f(x) dx,\quad {z\in[-\kappa,\kappa]}.
$$
The central issue is the effective computation of the eigenfunctions of the integral operator
\begin{equation}
\label{ee}
(EE^*f)(z)=2\int_{-\kappa}^\kappa\frac{\sin \tau(z-w)}{z-w}f(w) dw, \; {z\in[-\kappa,\kappa]}.
\end{equation}
This problem was beautifully solved by Landau, Pollak and Slepian  \cite{SP, LP1} in the early 1960's
by showing that the integral operator in \eqref{ee} commutes with the differential operator
$$
R(z, \partial_z)= \partial_z (\kappa^2-z^2)\partial_z -\tau^2,
$$ 
from which they described the common eigenfunctions via the differential operator. 
Note that $R(z,\partial_z)$ is the ``radial part'' of the Laplacian in prolate-spheroidal coordinates, motivating our title.
%Remarkably, this commuting property even appeared in the classical text by Ince \cite{Ince}. 
The commuting property was used by Fuchs \cite{Fuchs64} and Slepian \cite{Slepian65} to carry out a detailed analysis of the 
asymptotics of the eigenvalues of $E E^*$, while Jimbo et al. \cite{JMMS80} showed that 
its Fredholm determinant is a $\tau$-function of Painlev\'e V.

Remarkably, this commuting property appeared as early as 1907 in the work of Bateman \cite[Eqn. 38-41 accompanied by some differentiation]{bateman} and later in the classical text by Ince \cite{Ince}.
Mehta  \cite{Mehta} independently discovered and used it to analyze the Fredholm determinant of 
the integral operator \eqref{ee}, which he then applied to asymptotic problems in random matrices.
For recent numerical work on prolate spheroidal operators see \cite{osipov}; for applications to geophysics see \cite{simons,Fuchs64,Slepian65,JMMS80}.

Slepian \cite{S} found an extension of the time-band limiting analysis to $n$-dimensions. His method 
was based on passing to polar coordinates and then relying on a different commutativity result. 
He proved that for integer $N$ the integral operator
\[
(\mathcal Ef)(z) = \int_0^1 J_N(czw)\sqrt{czw}f(w)dw
\]
acting on a subspace of $L^2(0, 1; dw)$ with appropriate boundary conditions admits the commuting differential operator 
\[
\partial_z (1-z^2)\partial_z - c^2z^2+\frac{\frac{1}{4}-N^2}{z^2},
\]
where $J_N(x)$ denote the Bessel functions of the first kind.

In the early 1990's Tracy and Widom \cite{TW1,TW2} discovered 
one more remarkable commuting pair of integral and differential operators
associated to the Airy kernel. They effectively used this pair and a modification of the one for the Bessel kernel
in their study of the asymptotics of the level spacing distribution functions 
of the edge scaling limits of the Gaussian Unitary Ensemble and the 
Laguerre and Jacobi Ensembles. More precisely, Tracy and Widom 
proved that the integral operator with the Airy kernel
\[
\frac{A(z) A'(w) - A'(z) A(w)}{z-w}
\]
acting on $L^2(\tau, + \infty;dw)$
admits the commuting differential operator
\[
\partial_z (\tau-z) \partial_z - z(\tau-z),
\]
where $A(z)$ denotes the Airy function.

All of the above developments fit into one general scheme: 
commuting differential operators were constructed for an integral kernel of the form
\begin{equation}
\label{K-oper}
K_{\psi}(z,w) := \int_{\Gamma_2} \psi(x,z)\psi^*(x,w) dx
\end{equation}
acting on $L^2(\Gamma_1;dw)$, where $\Gamma_1$ and $\Gamma_2$ are contours in ${\mathbb{C}}$, $\psi(x,z)$ is a 
wave function for the KP hierarchy, $\psi^*(x,z)$ is its adjoint wave function.
Note in Slepian's Bessel-type example above, we get this kernel form for the square of his integral operator.
Many other instances of such commuting pairs were later discovered \cite{grunbaum-musings,Grunbaum-CPAM94,Grunbaum-AlgAnal96,CY18}, to name a few, and generalized to discrete and matrix-valued settings \cite{GPZ,Grunbaum-CMP2018}.

\subsection{The KP and KdV hierarchies} The Korteweg--de Vries (KdV) equation 
$$
\partial_t u  + u \partial_x u + \partial_x^3 u =0
$$
was introduced more than a century ago to model waves on shallow water surfaces. Its 
complete integrablity was established by Miura--Gardner--Kruskal \cite{MGK} and 
Lax \cite{Lax}. A wave function for a solution $u(x,t)$ is a function $\psi(x,z;t)$ satisfying
\[
(\partial_x^2 - u(x,t) )\psi(x,z;t) = z^2 \psi(x,z;t).
\]
The KdV equation fits into an infinite system of completely integrable nonlinear partial differential equations in variables $x,t_0,t_1,t_2,\dots$ known as the KP hierarchy.  Alternatively the KdV equation fits into the KdV hierarchy describing KP solutions independent of even times.

The KP hierarchy is an infinite dimensional integrable system whose wave functions $\psi(x,z)$ are eigenfunctions of differential operators 
$L(x, \partial_x)$ of higher order, 
and more generally of formal pseudo-differential operators. We refer the reader to van Moerbeke's exposition of the subject  \cite{vanMoerbeke} 
from the point of view of evolution on the (infinite-dimensional) Sato's Grassmannian ${\mathrm{Gr}}^{\mathrm{Sato}}$ 
and its applications to quantum gravity and intersection theory 
on moduli spaces of curves via the Kontsevich theorem \cite{Kontsevich}. The latter concerns 
precisely the solution of the KP hierarchy corresponding to the Airy wave function $\psi_{\mathrm{Ai}}(x,z) = A(x+z)$.

In the late 1970's Airault, McKean and Moser \cite{AMM-CPAM77} found a remarkable connection 
between the (infinite dimensional) KdV equation and finite dimensional 
integrable systems. They proved that any rational solution of the KdV equation that vanishes at infinity  has 
the form
\[
u(x, t) = \frac{1}{2}\sum_{i=1}^n \frac{1}{(x - x_i(t))^2}
\]
and that the KP flow for $t =t_1$ corresponds to the motion of the poles 
$(x_1(t), \ldots, x_n(t))$ according to the Calogero--Moser system with Hamiltonian
$H = \sum_i p_i^2/2 - \sum_{i<j} (x_i - x_j)^{-2}$. Krichever \cite{Krichever} proved
that this is true for every rational solution of KP vanishing at infinity
and that all solutions of the Calogero--Moser system arise in this way.

\subsection{Bispectrality and the Adelic Grassmannian}
The bispectral problem, posed by Duistermaat and the second named author in \cite{DG86}, 
asks for a classification of all functions $\psi(x,z)$ on a subdomain $\Omega_1 \times \Omega_2 \subseteq {\mathbb{C}^2}$ for which there exist 
two differential operators  $L(x, \partial_x)$ and $\Lambda(z, \partial_z)$ on 
$\Omega_1$ and $\Omega_2$ and two functions $\theta \colon \Omega_1 \to {\mathbb{C}}$, $f \colon \Omega_2 \to {\mathbb{C}}$,
such that
\begin{align*}
L(x, \partial_x) \psi(x,z) &= f(z) \psi(x,z),  
\\
\Lambda(z, \partial_z) \psi(x,z)
&= \theta(x) \psi(x, z).
\end{align*}
Many important relations of bispectrality to representation theory, algebraic and noncommutative geometry were 
subsequently found. %\cite{CRM-volume}. 
Early on it was realized that it is advantageous to think of its solutions as wave 
functions of the KP hierarchy. In this setting \cite{DG86} provided a classification of the 
second order bispectral operators $L(x,\partial_x)$.  Half of these come from rational solutions of the KdV equation.
The other half is comprised of the
Airy wave function $\psi_{\mathrm{Ai}}(x,z) = A(x+z)$,
the Bessel wave functions $\psi_{\mathrm{Be(\nu)}}(x,z) = \sqrt{xy} J_\nu(\sqrt{xy})$ for $\nu\in\mathbb{C}\backslash\mathbb{Z}$, and wave functions obtained from them by ``master symmetries'' of the KdV hierarchy \cite{ZM}.

Wilson made a deep insight to the bispectral problem \cite{Wilson-93}, providing the concept of classifying bispectral functions $\psi(x,z)$ according to their rank, defined as the greatest common divisor of the orders of all differential operators $L(x,\partial_x)$ having $\psi(x,z)$ as an eigenfunction.
For example, in the order 2 classification of \cite{DG86}, the wave functions of the rational solutions of KdV are of rank 1, 
while the remaining families are of rank 2.

In \cite{Wilson-93} Wilson classified all bispectral functions $\psi(x,z)$ of rank 1 in terms of an infinite dimensional sub-Grassmannian ${\mathrm{Gr}}^{\mathrm{ad}}$ of Sato's 
Grassamannian ${\mathrm{Gr}}^{\mathrm{Sato}}$, called the Adelic Grassmannian.
${\mathrm{Gr}}^{\mathrm{ad}}$ consists of those planes $W \in {\mathrm{Gr}}^{\mathrm{Sato}}$ obtained 
from the base plane $W_0 = {\mathbb{C}}[z]$ by imposing ``adelic-type'' conditions at finitely many points.
It was shown 
in \cite{BHY97} that these are precisely the KP wave functions $\psi(x,z)$ such that 
$$\psi(x,z) = \frac{1}{p(x)q(z)}P(x,\partial_x)\cdot e^{-xz}$$
and 
$$e^{-xz} = \frac{1}{\widetilde{p}(x) \widetilde{q}(y)}  \widetilde{P}(x,\partial_x) \cdot \psi(x,z)$$
for some differential operators $P(x,\partial_x)$ and $\widetilde{P}(x,\partial_x)$ with polynomial coefficients and polynomials
$p(x)$, $\widetilde{p}(x)$, $q(z)$, $\widetilde{q}(y)$. The orders of the differential operators 
$P$ and $\widetilde{P}$ will be called \vocab{degree} and \vocab{codegree} of $\psi(x,z)$, respectively.

Wilson \cite{Wilson-98} completed the circle back to Airault-McKean-Moser \cite{AMM-CPAM77} and Krichever \cite{Krichever}
by showing that the Adelic Grassmannian is the disjoint union of the Calogero-Moser spaces
${\mathrm{CM}}_n \subset {\mathrm{Gr}}^{\mathrm{Sato}}$ which are compactifications of the phase spaces of the 
Calogero-Moser integrable systems on the rational solutions of the KP hierarchy of \cite{AMM-CPAM77,Krichever}
\begin{equation}
\label{decomp}
{\mathrm{Gr}}^{\mathrm{ad}} = \bigsqcup_{n \geq 1}
{\mathrm{CM}}_n.
\end{equation}
%where each stratum ${\mathrm{CM}}_n$ is itself a obtained by a compactification 
%of the hamiltonian reduction picture for the Calogero-Moser phase space of Kazhdan-Kostant-Sternberg \cite{KKS}.

\section{Integral operators and points of ${\mathrm{Gr}}^{\mathrm{ad}}$}
\subsection{Reflectivity} 
The unifying feature of the diverse lines of research described above is a collection of hand-made examples of integral operators with kernels of the form in \eqref{K-oper} commuting with differential operators, obtained from certain specific wave functions $\psi(x,z)$ in ${\mathrm{Gr}}^{\mathrm{ad}}$.
%The most general class of examples known in previous literature were the $\psi(x,z)$ satisfying the two symmetry 
%constraints $\psi(x,z)=\psi(-x,-z)$ and $\psi^*(x,z) = \psi(x,z)$ \cite{CY18}.

For a long time, the examples provided above were the only known examples, and for this reason it was tempting to believe that it was a complete collection of examples.
However, this is not true at all!
In this paper we give a general solution of the problem that is applicable to the integral operators associated to the wave functions of \textbf{all} points of the Adelic Grassmannian.
It is based on a conceptual way of constructing the commuting differential operators from bispectral algebras.
Our key idea is that the intrinsic property of all of these integral operators is a more general one 
than a naive commutativity:
\begin{de}\label{Def-reflect}
An integral operator $T$, acting on $L^2(\Gamma)$ for a contour $\Gamma \subset {\mathbb{C}}$, 
is said to \vocab{reflect} a differential operator $R(z,\partial_z)$ if 
$$T\circ R(-z,-\partial_z) = R(z,\partial_z)\circ T$$
on a dense subspace of $L^2(\Gamma)$.
\end{de}
In the special case that a wave function $\psi(x,z) \in {\mathrm{Gr}}^{\mathrm{Sato}}$ satisfies the symmetry condition
$\psi(x,z) = \psi(-x,-z)$, this property for the kernel in \eqref{K-oper}
reduces to classical commutativity. 
This happens, for example, in the case of master symmetries \cite{Grunbaum-AlgAnal96}.
However, we will show that even more generally, imaginary rotation arguments 
transform reflecting pairs to classically commuting ones.

\begin{rem}
The reflection identity of Definition \ref{Def-reflect} is sensitive to the extension of the differential operator $R(z,\partial_z)$ to $L^2(\Gamma)$, which is not unique, and may hold for a unique choice of this extension.
This is a technical point that is often omitted in the classical prolate-spheroidal picture \cite{katsnelson}.
\end{rem}

%We expect that, for the purposes of spectral analysis of integral operators, 
%the reflectivity property may be exploited in a way similar to the case of classical commutativity.

\subsection{First general theorem -- reflection vs commutation}
Our first theorem associates to \emph{any} wave function $\psi(x,z)\in {\mathrm{Gr}}^{\mathrm{ad}}$ an integral operator $T$ which reflects a differential operator.
The reflected differential operator $R(z,\partial_z)$ resides in a natural algebra of differential operators associated to $\psi(x,z)$, called the \vocab{(right) generalized Fourier algebra}, defined in \cite{CY18} by
%$$\mathcal F_x(\psi) = \{L: \exists\ R\ \text{with}\ L(x,\partial_x)\psi(x,z) = R(z,\partial_z)\psi(x,z)\}.$$
$$\mathcal F_z(\psi) := \{R: \exists\ L\ \text{with}\ L(x,\partial_x)\psi(x,z) = R(z,\partial_z)\psi(x,z)\}.$$
The differential operators $L(x, \partial_x)$ that appear in the left hand side also form an algebra,
called  the \vocab{(left) generalized Fourier algebra} and denoted by $\mathcal F_x(\psi)$.
The map $L(x,\partial_x) \mapsto R(z, \partial_z)$ defines an algebra 
antiisomorphism 
$$
b_\psi :  \mathcal F_x(\psi) \rightarrow \mathcal F_z(\psi).
$$
The algebras $\mathcal F_x(\exp(-xz))$ and $\mathcal F_z(\exp(-xz))$ are both equal to the first Weyl algebra and the corresponding 
map $b$ is closely related to the Fourier transform.
%These algebras are naturally described by differential operators on line bundles over algebraic curves.

\begin{thm}\label{theorem 1}
For every wave function $\psi(x,z) \in {\mathrm{Gr}}^{\mathrm{ad}}$, the integral operator
$T_\psi$ on $L^2[t, \infty)$
%T: f(z)\mapsto \int_t^\infty K_\psi(z,w)f(w)dw
with kernel 
\begin{equation}
K_{\psi}(z,w) := \int_s^{\infty} \psi(y,z)\psi^*(y,w) dy
\label{K-oper-R}
\end{equation}
reflects a (non-constant) differential operator $R(z,\partial_z)\in\mathcal F_z(\psi)$
of order at most $2 \min(d_1, d_2)$ where $d_1$ and $d_2$ are the degree and codegree
of $\psi(x,z)$.
\end{thm}
A key feature of the proof of the theorem, sketched below, explicitly reduces the problem of finding the operator $R(z,\partial_z)$ to a finite-dimensional linear algebra problem.
This in turn provides an {\bf{effective algorithm}} for computing the reflected differential operator for all $\psi(x,z) \in {\mathrm{Gr}}^{\mathrm{ad}}$. 
In particular, we obtain examples of integral operators commuting with differential operators of orders much higher than can be reasonably found by hand, as shown in Examples \ref{first nontrivial example} and \ref{grunbaum operator example} below.

\subsection{Wilson's three involutions}
In general the operator $T$ defined by Theorem \ref{theorem 1} is not self-adjoint (even formally). In this way 
we may gain additional insight into the spectra of non-self adjoint integral operators.
In connection to Shannon's original questions, 
we have to be able to detect which operators in Theorem \ref{theorem 1} are of the form $EE^*$ 
and in particular are self-adjoint. For this we consider the three natural involutions of the Adelic Grassmannian ${\mathrm{Gr}}^{\mathrm{ad}}$ 
introduced by Wilson in \cite{Wilson-93}, along with a fourth involution not previously featured in this context corresponding to Schwartz reflection.
\begin{center}
\begin{tabular}{|c|c|}\hline
Name & Involution\\\hline
Adjoint & $a(\psi)(x,z) = \psi^*(x,z) = \frac{\wt P^*(x,\partial_x)\cdot e^{-xz}}{p(x)\wt p(x)}$\\\hline
Bispectral & $b(\psi)(x,z)=\psi(z,x)$\\\hline
Sign & $s(\psi)(x,z) = \psi(-x,-z)$\\\hline
Schwartz & $c(\psi)(x,z) = \ol{\psi(\ol x,\ol z)}$\\\hline
\end{tabular}
\end{center}
Note that the adjoint involution was used implicitly in \eqref{K-oper}. Wilson observed that the 
involutions $a$, $b$ and $s$ have the remarkable property that $a b$ is not an involution, but rather
\begin{equation}
\label{abs}
(ab)^2 =s.
\end{equation}

\subsection{Sketch of the proof of Theorem 1}\mbox{}

{\em{Step 1.}} Another way to phrase Wilson's property in \eqref{abs} is that
\begin{equation}
\label{Wilson-prop-2}
b_{a\psi}(b_\psi^{-1}(R)^*)^* (z,\partial_z) = R(-z,-\partial_z), \quad 
\forall R \in \mathcal F_z(\psi).
\end{equation}
Consider a differential operator $R_{s,t}(z,\partial_z)\in\mathcal F_z(\psi)$ such that both bilinear concomitants 
$$
\mathcal C_{b_\psi^{-1} R_{s,t}}(f,g;s) \quad \mbox{and} \quad \mathcal C_{R_{s,t}}(f,g;-t)
$$
are identically zero. We refer the reader to \cite{Wilson-concon} for the definition and properties of bilinear 
concomitants of differential operators.  Applying the identity in \eqref{Wilson-prop-2} together with integration by parts and the maps $b_\psi^{-1}$ and 
$b_{a \psi}$, we obtain that such an operator $R_{s,t}(z,\partial_z)$ satisfies
\begin{align*}
&R_{s,t}(z,\partial_z)\cdot K_\psi(z,w) = \\
  & = \int_s^\infty \big( R_{s,t}(z,\partial_z)\cdot\psi(x,z) \big) \psi(x,w)^* dx\\
  & = \int_s^\infty \big( b_\psi^{-1}(R_{s,t})(x,\partial_x)\cdot\psi(x,z) \big) \psi(x,w)^* dx\\
  & = \int_s^\infty \psi(x,z) \big( b_\psi^{-1}(R_{s,t})^*(x,\partial_x)\cdot\psi(x,w)^* \big) dx\\
  & = \int_s^\infty \psi(x,z) \big( b_{a\psi}(b_\psi^{-1}(R_{s,t}))^*(w,\partial_w)\cdot\psi(x,w)^* \big) dx\\
  & = R_{s,t}^*(-w,-\partial_w)\cdot K_\psi(z,w).
\end{align*}
This identity combined with one more integration by parts proves that
$$R_{s,t}(z,\partial_z) \circ T_\psi  = T_\psi \circ R_{s,t}(-z,-\partial_z).$$
for the integral operator $T_\psi$ with kernel as in \eqref{K-oper-R}. 

The remainder of the proof of Theorem $1$ revolves around demonstrating the existence of a differential operator $R_{s,t}(z, \partial_z) \in \mathcal F_z(\psi)$ 
satisfying the conditions of Step 1.
Its existence, along with a sharp upper bound on its order, are obtained 
by algebro-geometric arguments. 

{\em{Step 2.}} The operators in the Fourier algebra $\mathcal F_x(\psi)$ 
naturally have a co-order ${\mathrm{coord}} R(z, \partial_z) := {\mathrm{ord}} (b_\psi^{-1} R)(x, \partial_x)$.
For a pair of nonnegative integers $\ell, m$, set
$$
\mathcal F_z^{\ell,m}(\psi) := \{ R \in \mathcal F_z(\psi) : {\mathrm{ord}} R \leq \ell, {\mathrm{coord}} R \leq m \}.
$$
Recall the decomposition in \eqref{decomp}; let $\psi(x,z) \in {\mathrm{CM}}_n \subset {\mathrm{Gr}}^{\mathrm{ad}}$.
One shows that $\mathcal F_z (\psi)$ 
is isomorphic to the algebra of differential operators on a rank 1, torsion free sheaf over the spectral curve of the solution of KP with 
wave function $\psi(x,z)$. Interpreting $n$ as the \emph{differential genus} of the sheaf of the curve in the sense of Berest-Wilson \cite{BW04} 
and then converting it to the Letzter-Markar-Limanov invariant of the sheaf shows that
\begin{align*}
\dim \mathcal F_z^{\ell, m} (\psi) &= (\ell +1) (m +1) - n 
\\
&\geq (\ell+1)(m+1) - 2 \min(d_1, d_2)^2
\end{align*}
for $\ell, m \geq  2 \min(d_1, d_2) -1$.

{\em{Step 3.}} For a differential operator $R(z, \partial_z)$ of order $\leq \ell$, the 
identical vanishing of the concomitant $C_R(f,g;-t)$ is shown to lead to at most 
$\lceil \ell/2 \rceil . \lceil (\ell+1)/2 \rceil$ linearly independent (linear) conditions on 
the coefficients of $R$ and their derivatives. 

This estimate, combined with that in Step 2, 
proves the existence of a a differential operator $R_{s,t}(z, \partial_z) \in \mathcal F_z(\psi)$ 
satisfying the conditions of Step 1 of order at most $2 \min (d_1, d_2)$.
\qed
%\begin{align*}
%R(z,\partial_z)\cdot T(f(w))(z)
%  & = \int_t^\infty R(z,\partial_z)\cdot K(z,w)f(w)dw\\
%  & = \int_t^\infty R(-w,-\partial_w)^*\cdot K(z,w)f(w)dw\\
%  & = \int_t^\infty K(z,w)R(-w,-\partial_w)\cdot f(w)dw\\
%  & = T(R(-w,-\partial_w)\cdot f(w))(z).
%\end{align*}
\begin{rem} {\em{(i)}} Wilson's identity on involutions in the Adelic Grassmannian in \eqref{abs} and its use 
in Step 1 are the intrinsic reasons for the appearance of reflectivity in Theorem \ref{theorem 1} 
rather than classical commutativity.

{\em{(ii)}} All previous approaches for constructing commuting pairs of integral and differential 
operators, like those in \cite{SP, LP1, S,TW1,TW2}, relied on a by-hand construction of a commuting 
differential operator. Step 1 of the proof is where bispectrality plays a deep role and the operator 
is constructed from the generalized Fourier algebra $\mathcal F_z(\psi)$. 
\end{rem}

\section{The Laplace vs Fourier pictures}
\subsection{Second general theorem -- the Laplace picture}
Consider a wave function $\psi \in {\mathrm{Gr}}^{\mathrm{ad}}$.
We draw a parallel between the integral operators from Theorem \ref{theorem 1} and those of the form $EE^*$ by considering the following analogs of the Laplace transform and its adjoint:
\begin{align*}
L_\psi &: f(x)\mapsto \int_0^\infty \psi(y,z)f(y)dy,\\
L^*_\psi&: f(z)\mapsto \int_0^\infty \ol{\psi(x,w)}f(w)dw.
\end{align*}
In the special case that $\psi(x,z) = \exp(-xz)$, the operator $L_\psi$ is precisely the Laplace transform.
The \vocab{time and band-limited} versions of these are (for $z\geq t$)
$$({\mathcal{E}}_\psi f)(z)   = (\chi_{[t,\infty)}L_\psi \chi_{[s,\infty)} f)(z)  = \int_s^\infty \psi(y,z) f(y) dy,$$
and (for $x\geq s$)
$$({\mathcal{E}}^*_\psi f)(x) = (\chi_{[s,\infty)}L^*_\psi \chi_{[t,\infty)} f)(x) = \int_t^\infty \ol{\psi(x,w)}f(w)dw.$$
They give rise to the self-adjoint operator analogous to the one considered by Landau, Pollak, and Slepian
\begin{align*}
({\mathcal{E}}_\psi {\mathcal{E}}^*_\psi f)(z) &= \int_t^\infty K_\psi(z,w)f(w)dw, \; \; \mbox{where} \\
K_\psi (z,w) &= \int_s^\infty \psi(y,z)\ol{\psi(y,w)}dy,
\end{align*}
viewed as an operator on $L^2(t, \infty)$. 
Under natural mild conditions on $\psi(x,z)$, Theorem \ref{theorem 1} determines the existence of differential operators reflected by ${\mathcal{E}}_\psi {\mathcal{E}}^*_\psi$.
For a different situation involving the Laplace transform, see \cite{BG85}.
\begin{thm}\label{theorem 2}
For every wave function $\psi(x,z)$ in Wilson's adelic Grassmannian, fixed under the involution $ac$ of ${\mathrm{Gr}}^{\mathrm{ad}}$ (defined by the table of involutions above),  
the integral operator ${\mathcal{E}}_\psi {\mathcal{E}}^*_\psi$ reflects a (non-constant) 
differential operator $R(z,\partial_z)\in\mathcal F_z(\psi)$ 
of order at most $2 \min (d_1, d_2)$ 
where $d_1$ and $d_2$ are the degree and codegree of $\psi(x,z)$. 
\end{thm}
\noindent
\begin{proof}[Sketch of Proof] From the assumption that $\psi(x,z)$ is fixed under 
the involution $ac$, one deduces that $\psi^*(x,z) = \overline{\psi(x,z)}$ 
for $x,z \in {\mathbb{R}}$. From this one shows that ${\mathcal{E}}_\psi {\mathcal{E}}^*_\psi$ 
equals the integral operator with kernel $K_\psi$ from Theorem \ref{theorem 1}.
\end{proof}
\begin{rem}
Under the assumption that $\psi(x,z)$ is fixed by $ac$, the reflected operator $R_{s,t}(z,\partial_z)$ satisfies the identity $R_{s,t}^*(z,\partial_z) = R_{s,t}(-z,-\partial_z)$.
In this case, the reflection property may be restated in the form
$${\mathcal{E}}_\psi {\mathcal{E}}^*_\psi \circ R_{s,t}^*(z,\partial_z) = R_{s,t}(z,\partial_z)\circ {\mathcal{E}}_\psi {\mathcal{E}}^*_\psi.$$
\end{rem}
\begin{ex}\label{easy ex} 
Consider the simplest case $\psi(x,z) = \exp(-xz)$. 
The integral operator
\begin{align*}
({\mathcal{E}}_\psi {\mathcal{E}}^*_\psi f)(z) = \int_t^\infty \frac{\sinh(s(z+w))}{z+w}f(w)dw\\
\end{align*}
acting on $L^2(t, \infty)$
reflects the first order differential operator 
$$R_{s,t}(z,\partial_z) = (z + t)\partial_z + sz.$$ 
All previous works on this kernel deal with a commuting second order
differential operator.
\end{ex}
%A more complicated example not featured before in the literate is the following.
The wave functions associated to rational solutions \cite{AMM-CPAM77} of KdV are automatically fixed by the involution $a$.
Additionally, those with real coefficients are fixed by $c$ and thus satisfy the assumptions 
of Theorem \ref{theorem 2}. These are precisely the bispectral functions 
in the KdV family in \cite{DG86} with real coefficients
(associated to second order differential operators of rank 1).
There has been a substantial effort since 1986 to find commuting 
differential operators for the corresponding integral operators, but absolutely 
no examples have been found beyond the case $\psi(x,z) = \exp(- xz)$ or \cite{Grunbaum-AlgAnal96}.
The next example demonstrates how Theorem \ref{theorem 2} 
resolves this problem. 
\begin{ex}\label{first nontrivial example}
Let $r \in {\mathbb{R}}^*$. Consider the function
$$\psi(x,z) = \frac{(x+z^{-1})^3 - z^3 + r}{x^3 + r}e^{-xz},$$
which up to a change of variables is precisely the first 
nontrivial bispectral function in \cite{DG86} given on Eq. (1.39).
The integral operator ${\mathcal{E}}_\psi {\mathcal{E}}^*_\psi$ has kernel
$$K_\psi(z,w) = \frac{\psi(s,z)\psi_x(s,w)-\psi_x(s,z)\psi(s,w)}{z^2-w^2}.$$
By Theorem \ref{theorem 2} it reflects a differential operator in $\mathcal F_z(\psi)$.
Our algorithm produces an operator of order $3$, given by
\begin{align*}
R_{s,t}(z,\partial_z)
  &= -(z+t)^2z\partial_z^3 + \big(st^3 - 3stz^2 - 2sz^3 - t^3rz^2\\
  & - t^2rz^3 - \frac{3}{2}t^2 - 6tz - \frac{9}{2}z^2 \big)\partial_z^2 + \big(s^2t^2z - s^2z^3\\
  & - 2st^3rz^2 - 2st^2rz^3 - 6stz - 6sz^2 - 2t^3rz\\
  & - 3t^2rz^2 + 6t^2z^{-1} + 6t - z \big)\partial_z -6st^3z^{-2} - 3t^2z^{-2}\\
  & + s^3tz^2 - s^2t^3rz^2 - s^2t^2rz^3 - \frac{3}{2}s^2z^2 - 2st^3rz\\
  & - 3st^2rz^2 - t^2rz + 3.
\end{align*}
\end{ex}

\subsection{Third general theorem -- the Fourier picture}
By performing a \vocab{90 degree rotation} in the complex variable $z$, we move from the Laplace transform picture to the Fourier transform picture.
We prove that in this way one can {\em{convert the reflected differential operators in Laplace picture to commuting differential operators 
in the Fourier picture}}.
Specifically, we replace the operators $L_\psi$ and $L_\psi^*$ with their Fourier counterparts
\begin{align*}
F_\psi &: f(x)\mapsto \int_{-\infty}^\infty \psi(y,-iz)f(y)dy,\\
F_\psi^*&: f(z)\mapsto \int_{-\infty}^\infty \ol{\psi(x,-iw)}f(w)dw.
\end{align*}
In the special case $\psi(y,z) = \exp(-yz)$, the operator $F_\psi$ is the Fourier transform.
We define the time and band-limited operators $E_\psi$ and $E^*_\psi$ similarly to ${\mathcal{E}}_\psi$ and ${\mathcal{E}}^*_\psi$:
\begin{align*}
&(E_\psi f)(z)  = (\chi_{[t,\infty)} F_\psi \chi_{[s,\infty)} f)(z)  = \int_s^\infty \psi(y,-iz) f(y) dy, \\
&(E^*_\psi f)(x) = (\chi_{[s,\infty)}F^*_\psi \chi_{[t,\infty)} f)(x) = \int_t^\infty \ol{\psi(x,-iw)}f(w)dw.
\end{align*}
The self-adjoint operator
$$(E_\psi E^*_\psi f)(z) = \int_s^\infty \int_t^\infty \psi(y,-iz)\ol{\psi(y,-iw)}f(w)dwdy$$
acting on $L^2(t, \infty)$ 
no longer has a simple kernel expression as above since the relevant integral does not converge outright, but can be given sense as a distribution.
Even so, the method of proof of Theorem \ref{theorem 1} applies, giving us a certain relationship between an integral and differential operators.
Serendipitously, due to the change in sign with complex conjugation, in this case we obtain a strict commutativity relation.
\begin{thm}\label{theorem 3}
For every wave function $\psi(x,z)$ in Wilson's adelic Grassmannian, fixed under the involution $ac$ of ${\mathrm{Gr}}^{\mathrm{ad}}$, 
the integral operator $E_\psi E_\psi^*$ commutes with the differential operator
$R_{s, it}(-iz, i \partial_z)$ where $R_{s,t}(z, \partial_z)$ is the corresponding differential operator in Theorem \ref{theorem 2} 
(its coefficients are rational functions in $z, s, t$). 
\end{thm}
In particular, we obtain that $E_\psi E^*_\psi$ commutes with the self-adjoint operator 
$R_{s, it}(-iz, i \partial_z) R_{s, it}^*(-iz, i \partial_z)$.
\begin{proof}[Sketch of Proof] One repeats Step 1 of the proof of Theorem \ref{theorem 1}
to show that $R_{s,it}( - i z, i \partial_z)$ commutes with  $E_\psi E_\psi^*$ 
for every differential operator $R(z,\partial_z)\in\mathcal F_z(\psi)$ 
for which both bilinear concomitants 
$$
\mathcal C_{b_\psi^{-1} R_{s,it}(-iz,i\partial_z)}(f,g;s) \quad \mbox{and} \quad \mathcal C_{R_{s,it}(-iz,i\partial_z)}(f,g;it)
$$
are identically zero. The operator $R_{s,t}(z, \partial_z)$ from Theorem \ref{theorem 2}
has these properties because its coefficients are rational functions in $z, s, t$. 
This is proved by an analysis of the structure of the algebra $\mathcal F_z(\psi)$.
\end{proof}

Note that for analytic reasons one cannot deduce Theorem \ref{theorem 3} from 
Theorem \ref{theorem 2} by an elementary change of variables.

Theorems \ref{theorem 1}--\ref{theorem 3} form the foundation for our forthcoming series of papers on the asymptotics of the eigenvalues of the integral operators associated to the wave functions of the rational solutions of the KP equation and the numerical properties of the associated eigenfunctions (which generalize the prolate spheroidal wave functions).

\begin{ex}
Consider the case $\psi(x,z)=\exp(-xz)$ as in Example \ref{easy ex}.
The self-adjoint integral operator is given by
$$(E_\psi E^*_\psi f)(z) = \int_s^\infty \int_t^\infty e^{iy(z-w)}f(w)dwdy.$$
{\em{The eigenvalues of this operator are precisely the 
singular values of the semi-infinite time-band limiting of the Fourier transform.}}
This integral operator commutes with the first order differential operator 
$$R_{s,it}(-iz,i \partial_z) = (z - t)\partial_z - isz$$
obtained from the differential operator in Example \ref{easy ex}.
As a consequence we obtain that $E_\psi E^*_\psi$ commutes with the self-adjoint second order differential operator
$$-R_{s,it}(-iz,i\partial_z) R_{s,it}^*(-iz,i\partial_z) = \partial_z(z-t)^2\partial_z - is\{z(z-t),\partial_z\},$$
where here $\{\cdot,\cdot\}$ denotes the anti-commutator bracket.
\end{ex}

%Although time-band limiting in the settings of the Fourier and Laplace transforms has been studied in great detail after  \cite{SP, LP1}, 
%we know of no previous observations (in the case of semi-infinite truncations) of commuting
%second order self-adjoint differential operators of the form $R R^*$.

%Remarkably, 
%a similar phenomenon holds 
%for all examples that we are aware of: 
%for a wave function $\psi \in {\mathrm{Gr}}^{\mathrm{ad}}$,
%among the self-adjoint differential operators of minimal order commuting with 
%$E_\psi E_\psi^*$ there is always one of the form $R R^*$, 
%where $R$ also commutes with $E_\psi E_\psi^*$.

\begin{ex}
Consider the wave function
$$\psi(x,z) = \frac{(x+z^{-1})^3 - z^3 + r}{x^3 + r}e^{-xz}$$
from Example \ref{first nontrivial example}.
The associated integral operator $E_\psi E^*_\psi$ commutes with the third order differential operator $R:=R_{s,it}(-iz,i \partial_z)$ 
where $R_{s,t}(z, \partial_z)$ is the differential operator in Example \ref{first nontrivial example}.
Its formal adjoint is $R^* = -R + s^2t^2 + 4$,
so that $E_\psi E^*_\psi$ commutes with the sixth order self-adjoint operator
$$-R R^* = R^2 - (s^2t^2 + 4)R.$$
%In fact, we prove in a further example below that it commutes with a self-adjoint fourth order differential operator.
\end{ex}

\section{Simultaneous reflectivity and commutativity}
The proof of Theorem \ref{theorem 1} produces a large algebra of reflected operators rather than a single one, 
because the argument can be applied to the full Fourier algebra $F_z(\psi)$ of $\psi \in {\mathrm{Gr}}^{\mathrm{ad}}$.
This can be used to prove the existence of \vocab{universal operators} which are simultaneously reflected by 
(or commute with) finite-dimensional collections of integral operators.
\begin{thm}\label{theorem 4}
{\em{(i)}} Consider any finite collection of wave functions $\{ \psi_k(x,z) : 1 \leq k \leq n \} \in  {\mathrm{Gr}}^{\mathrm{ad}}$ 
and let $T_k$ be the associated integral operators as in Theorem \ref{theorem 1} for the same values of $s$ and $t$.
There exists a non-constant differential operator in $\bigcap_k\mathcal F_z(\psi_k)$ simultaneously reflected by each of the integral operators $T_k$ for all $k$.

{\em{(ii)}} If, in addition, all wave functions $\psi_k(x,z)$ are fixed under the involution $ac$ of ${\mathrm{Gr}}^{\mathrm{ad}}$, then there exists 
a differential operator $R_{s,t}^{\text{univ}} (z, \partial_z)$ which is simultaneously reflected by all integral operators spanned by ${\mathcal{E}}_j {\mathcal{E}}_k^*$
for $1 \leq j, k \leq n$. This differential operator has rational coefficients in $z,s,t$ and $\wt R_{s,t}(z,\partial_z) := R_{s, it}^{\text{univ}}(-iz, i \partial_z)$ commutes with all
integral operators $E_j E_k^*$ for $ 1 \leq j, k \leq n$.
\end{thm}
In the situation of part (ii) all integral operators $E_j E_k^*$, $ 1 \leq j, k \leq n$, 
commute with the self-adjoint operator 
$\wt R_{s, t}(z, \partial_z) \wt R_{s, t}^*(z, \partial_z)$.
Furthermore, since the Fourier algebra of $\exp(-xz)$ is just the algebra of differential operators with polynomial coefficients, we can force all of the coefficients of 
$\wt R_{s,t}(z,\partial_z)$ to be polynomials in $z$.

\begin{ex}\label{grunbaum operator example}
Consider the pair of wave functions $\{\psi_1(x,z),\psi_2(x,z)\}$ with $\psi_n(x,z) = K_n(xz)\sqrt{xz}$ for $K_n(z)$ the modified Bessel function of the second kind.
Thus by Theorem \ref{theorem 4} there should exist 
a self-adjoint differential operator $\wt R_{s,t}(z,\partial_z)$ in $\mathcal F_z(\psi_1)\cap \mathcal F_z(\psi_2)$ with polynomial coefficients which commutes with the integral operators $E_kE_j^*$ defined by the wave functions $\psi_k(x,z)$ for $k=1,2$.
Note also that $\wt R_{s,t}(z, \partial_z)$ will commute with the integral operator $EE^*$ associated with the wave function from Example \ref{first nontrivial example} \textbf{for any $r$}, since this operator will be a linear combination of the $E_kE_j^*$'s. Using our algorithm for Theorem \ref{theorem 1}, we obtain an operator of order $6$ of the form
$$\wt R_{s,t}(z,\partial_z) = \sum_{m=0}^3 \partial_z^m f_m(z) \partial_z^m,$$
where
\begin{align*}
f_0(z) &= \frac{z^2(3s^6t^3-54s^4t)}{6} + s^6z^5 - \frac{3s^6tz^4}{2} + 12s^4z^3,\\
f_1(z) &= (z-t)\left(3s^4z^4-3s^4tz^3 + 12 s^2z^2 + 9s^2tz - 9s^2t^2\right),\\
f_2(z) &= (z-t)^2\left(3s^2z^3 - \frac{3s^2tz^2}{2} + 12t\right),\\
f_3(z) &= (z-t)^3z^2.
\end{align*}
\end{ex}

\section{Concluding remarks}
We have presented a unified general construction of commuting pairs based on the intrinsic properties of symmetries of soliton equations.
It has not escaped our notice that the specific connection we have described between commuting integral
and differential operators and solutions of the KdV equation, in
particular the critical role of the reflecting property in these
classical problems, opens up avenues of broad new applications of
integrable systems to spectral analysis of integral operators, going
far beyond sinc, Bessel and Airy kernels.
Additionally, the new pairs of commuting integral and differential operators may have a role to play in random matrix theory.

\subsection{Acknowledgments}
The research of I.Z. was supported by CONICET grant PIP 112-200801-01533 and SeCyT-UNC.
Additionally I.Z. and F.A.G. benefitted from conversations during the Research in Pairs program at Oberwolfach in the fall of 2018.
The research of M.Y. was supported by NSF grant DMS-1901830 and
Bulgarian Science Fund grant DN02/05.
The research of W.R.C. was supported by a 2018 AMS-Simons Travel Grant.

\bibliographystyle{abbrv}
\bibliographystyle{plain}
\bibliography{reflect}

\begin{thebibliography}{10}

\bibitem{AMM-CPAM77}
H.~Airault, H.~P. McKean, and J.~Moser.
\newblock Rational and elliptic solutions of the {K}orteweg-de {V}ries equation
  and a related many-body problem.
\newblock {\em Comm. Pure Appl. Math.}, 30(1):95--148, 1977.

\bibitem{BHY97}
B.~Bakalov, E.~Horozov, and M.~Yakimov.
\newblock Bispectral algebras of commuting ordinary differential operators.
\newblock {\em Comm. Math. Phys.}, 190(2):331--373, 1997.

\bibitem{bateman}
H.~Bateman.
\newblock On the inversion of a definite integral.
\newblock {\em Proceedings of the London Mathematical Society}, 2(1):461--498,
  1907.

\bibitem{BW04}
Y.~Berest and G.~Wilson.
\newblock Differential isomorphism and equivalence of algebraic varieties.
\newblock In {\em Topology, geometry and quantum field theory}, volume 308 of
  {\em London Math. Soc. Lecture Note Ser.}, pages 98--126. Cambridge Univ.
  Press, Cambridge, 2004.

\bibitem{BG85}
M.~Bertero and F.~A. Gr\"{u}nbaum.
\newblock Commuting differential operators for the finite {L}aplace transform.
\newblock {\em Inverse Problems}, 1(3):181--192, 1985.

\bibitem{CY18}
W.~R. Casper and M.~T. {Yakimov}.
\newblock {Integral operators, bispectrality and growth of Fourier algebras}.
\newblock {\em J. Reine Angew. Math}, 10.1515/crelle-2019-0031, 2019.

\bibitem{DG86}
J.~J. Duistermaat and F.~A. Gr{\"{u}}nbaum.
\newblock {Differential equations in the spectral parameter}.
\newblock {\em Comm. Math. Phys.}, 103(2):177--240, 1986.

\bibitem{Fuchs64}
W.~H.~J. {Fuchs}.
\newblock {On the eigenvalues of an integral equation arising in the theory of
  band-limited signals.}
\newblock {\em {J. Math. Anal. Appl.}}, 9:317--330, 1964.

\bibitem{Grunbaum-CPAM94}
F.~A. Gr\"{u}nbaum.
\newblock Time-band limiting and the bispectral problem.
\newblock {\em Comm. Pure Appl. Math.}, 47(3):307--328, 1994.

\bibitem{Grunbaum-AlgAnal96}
F.~A. Gr\"{u}nbaum.
\newblock Band-time-band limiting integral operators and commuting differential
  operators.
\newblock {\em Algebra i Analiz}, 8(1):122--126, 1996.

\bibitem{grunbaum-musings}
F.~A. Gr\"{u}nbaum.
\newblock Some bispectral musings.
\newblock In J.~Harnad and A.~Kasman, editors, {\em The bispectral problem
  ({M}ontreal, {PQ}, 1997)}, volume~14 of {\em CRM Proc. Lecture Notes}, pages
  31--45. Amer. Math. Soc., Providence, RI, 1998.

\bibitem{GPZ}
F.~A. Gr\"{u}nbaum, I.~Pacharoni, and I.~Zurri\'{a}n.
\newblock Bispectrality and time-band-limiting: Matrix valued polynomials.
\newblock {\em International Mathematics Research Notices},
  10.1093/imrn/rny140, 6 2018.

\bibitem{Grunbaum-CMP2018}
F.~A. Gr\"{u}nbaum, L.~Vinet, and A.~Zhedanov.
\newblock Algebraic {H}eun operator and band-time limiting.
\newblock {\em Comm. Math. Phys.}, 364(3):1041--1068, 2018.

\bibitem{Ince}
E.~L. Ince.
\newblock {\em Ordinary {D}ifferential {E}quations}.
\newblock Dover Publications, New York, 1944.

\bibitem{JMMS80}
M.~{Jimbo}, T.~{Miwa}, Y.~{M\^ori}, and M.~{Sato}.
\newblock {Density matrix of an impenetrable Bose gas and the fifth Painlev\'e
  transcendent.}
\newblock {\em {Physica D}}, 1(1):80--158, 1980.

\bibitem{katsnelson}
V.~Katsnelson.
\newblock Self-adjoint boundary conditions for the prolate spheroid
  differential operator.
\newblock In D.~Alpay and B.~Kirstein, editors, {\em Indefinite Inner Product
  Spaces, Schur Analysis, and Differential Equations}, pages 357--386.
  Springer, 2018.

\bibitem{Kontsevich}
M.~Kontsevich.
\newblock Intersection theory on the moduli space of curves and the matrix
  {A}iry function.
\newblock {\em Comm. Math. Phys.}, 147(1):1--23, 1992.

\bibitem{Krichever}
I.~M. Krichever.
\newblock Rational solutions of the {K}adomcev-{P}etvia\v{s}vili equation and
  the integrable systems of {$N$} particles on a line.
\newblock {\em Funkcional. Anal. i Prilo\v{z}en.}, 12(1):76--78, 1978.

\bibitem{LP1}
H.~J. Landau and H.~O. Pollak.
\newblock Prolate spheroidal wave functions, {F}ourier analysis and
  uncertainty. {II}.
\newblock {\em Bell System Tech. J.}, 40:65--84, 1961.

\bibitem{Lax}
P.~D. Lax.
\newblock Integrals of nonlinear equations of evolution and solitary waves.
\newblock {\em Comm. Pure Appl. Math.}, 21:467--490, 1968.

\bibitem{Mehta}
M.~L. Mehta.
\newblock {\em Random matrices}, volume 142.
\newblock Elsevier, 3 edition, 2004.

\bibitem{MGK}
R.~M. Miura, C.~S. Gardner, and M.~D. Kruskal.
\newblock Korteweg-de {V}ries equation and generalizations. {II}. {E}xistence
  of conservation laws and constants of motion.
\newblock {\em J. Mathematical Phys.}, 9:1204--1209, 1968.

\bibitem{osipov}
A.~Osipov, V.~Rokhlin, and H.~Xiao.
\newblock Prolate spheroidal wave functions of order zero.
\newblock {\em Springer Ser. Appl. Math. Sci}, 187, 2013.

\bibitem{Shannon1}
C.~E. Shannon.
\newblock A mathematical theory of communication.
\newblock {\em Bell system technical journal}, 27(3):379--423, 1948.

\bibitem{Shannon2}
C.~E. Shannon.
\newblock A mathematical theory of communication.
\newblock {\em Bell system technical journal}, 27(4):623--666, 1948.

\bibitem{simons}
F.~J. Simons, F.~Dahlen, and M.~A. Wieczorek.
\newblock Spatiospectral concentration on a sphere.
\newblock {\em SIAM review}, 48(3):504--536, 2006.

\bibitem{S}
D.~Slepian.
\newblock Prolate spheroidal wave functions, {F}ourier analysis and
  uncertainity. {IV}. {E}xtensions to many dimensions; generalized prolate
  spheroidal functions.
\newblock {\em Bell System Tech. J.}, 43:3009--3057, 1964.

\bibitem{Slepian65}
D.~{Slepian}.
\newblock {Some asymptotic expansions for prolate spheroidal wave functions.}
\newblock {\em {J. Math. Phys., Mass. Inst. Techn.}}, 44:99--140, 1965.

\bibitem{SP}
D.~Slepian and H.~O. Pollak.
\newblock Prolate spheroidal wave functions, {F}ourier analysis and
  uncertainty. {I}.
\newblock {\em Bell System Tech. J.}, 40:43--63, 1961.

\bibitem{TW1}
C.~A. Tracy and H.~Widom.
\newblock Fredholm determinants, differential equations and matrix models.
\newblock {\em Comm. Math. Phys.}, 163(1):33--72, 1994.

\bibitem{TW2}
C.~A. Tracy and H.~Widom.
\newblock Level-spacing distributions and the {A}iry kernel.
\newblock {\em Comm. Math. Phys.}, 159(1):151--174, 1994.

\bibitem{vanMoerbeke}
P.~van Moerbeke.
\newblock Integrable foundations of string theory.
\newblock In O.~Babelon, P.~Cartier, and Y.~Kosmann-Schwarzbach, editors, {\em
  Lectures on integrable systems ({S}ophia-{A}ntipolis, 1991)}, pages 163--267.
  World Sci. Publ., River Edge, NJ, 1994.

\bibitem{Wilson-concon}
G.~Wilson.
\newblock {On the antiplectic pair connected with the Adler-Gel'fand-Dikii
  bracket}.
\newblock {\em {Nonlinearity}}, 5(1):109--131, 1992.

\bibitem{Wilson-93}
G.~Wilson.
\newblock Bispectral commutative ordinary differential operators.
\newblock {\em J. Reine Angew. Math.}, 442:177--204, 1993.

\bibitem{Wilson-98}
G.~Wilson.
\newblock Collisions of {C}alogero-{M}oser particles and an adelic
  {G}rassmannian.
\newblock {\em Invent. Math.}, 133(1):1--41, 1998.
\newblock With an appendix by I. G. Macdonald.

\bibitem{ZM}
J.~P. Zubelli and F.~Magri.
\newblock Differential equations in the spectral parameter, {D}arboux
  transformations and a hierarchy of master symmetries for {K}d{V}.
\newblock {\em Comm. Math. Phys.}, 141(2):329--351, 1991.

\end{thebibliography}

\end{document}